\documentclass[aop,MSNbibl,nameyear,seceqn,dvips]{arximspdf}

\doi{10.1214/13-AOP850} 
\volume{42}
\issue{2}
\pubyear{2014}
\firstpage{794}
\lastpage{817}

\makeatletter
\newcommand{\rrVert}{\Vert}
\newcommand{\rrvert}{\vert}
\newcommand{\llVert}{\Vert}
\newcommand{\llvert}{\vert}
\def\cal{\mathcal}
\newtheorem{lemma}{Lemma}[section]

\newtheorem{theorem}{Theorem}[section]
\newtheorem{corollary}{Corollary}[section]
\newproclaim{remark}{Remark}[section]
\newproclaim{example}{Example}[section]
\newcommand{\etav}{\bolds{\eta}}
\date{}

\def\open{}
\def\pr{\mathsf{P}} %
\def\ep{\mathsf{E}} %
\def\N{{\open\mathbb{N}}}
\def\R{{\open\mathbb{R}}}
\def\Z{{\open\mathbb{Z}}}
\def\B{\mathbb{B}}
\def\eqd{\stackrel{\cal D}{=}}
\def\qed{\hfill$\diamondsuit$}
\makeatother

\begin{document}
\begin{frontmatter}

\title{Koml\'os--Major--Tusn\'ady approximation under dependence}
\runtitle{Koml\'os--Major--Tusn\'ady approximation}

\begin{aug}
\author[A]{\fnms{Istv\'an} \snm{Berkes}\corref{}\thanksref{t1}\ead[label=e1]{berkes@tugraz.at}},
\author[B]{\fnms{Weidong} \snm{Liu}\thanksref{t2}\ead[label=e2]{liuweidong99@gmail.com}}
\and
\author[C]{\fnms{Wei Biao} \snm{Wu}\thanksref{t3}\ead[label=e3]{wbwu@galton.uchicago.edu}}
\runauthor{I. Berkes, W. Liu and W. B. Wu}
\affiliation{Graz University of Technology, Shanghai Jiao Tong
University and University~of~Chicago}
\thankstext{t1}{Supported by FWF Grant P 24302-N18 and OTKA Grant
K~108615.}
\thankstext{t2}{Supported by NSFC Grant 11201298.}
\thankstext{t3}{Supported in part from DMS-09-06073 and
DMS-11-06970.}
\address[A]{I. Berkes\\
Institute of Statistics\\
Graz University of Technology\\
Kopernikusgasse 24\\
8010 Graz\\
Austria\\
\printead{e1}} 
\address[B]{W. Liu\\
Department of Mathematics\\
Shanghai Jiao Tong University\\
800 Dongchuan Road Minhang\\
Shanghai\\
China\\
\printead{e2}}
\address[C]{W. B. Wu\\
Department of Statistics\\
University of Chicago\\
5734 S. University Avenue\\
Chicago, Illinois 60637\\
USA\\
\printead{e3}}
\end{aug}

\received{\smonth{2} \syear{2012}}
\revised{\smonth{2} \syear{2013}}

%
\begin{abstract}
The celebrated results of Koml\'os, Major and Tusn\'ady
[\textit{Z. Wahrsch. Verw. Gebiete} \textbf{32} (1975) 111--131;
\textit{Z. Wahrsch. Verw. Gebiete} \textbf{34} (1976) 33--58]
give optimal Wiener approximation for the partial sums of i.i.d. random
variables and provide a powerful tool in probability and
statistics. In this paper we extend KMT approximation for a large
class of dependent stationary processes, solving a long standing
open problem in probability theory. Under the framework of
stationary causal processes and functional dependence measures of
Wu [\textit{Proc. Natl. Acad. Sci. USA} \textbf{102} (2005)
14150--14154],
 we show that, under natural moment conditions, the partial
sum processes can be approximated by Wiener process with an optimal
rate. Our dependence conditions are mild and easily verifiable. The
results are applied to ergodic sums, as well as to nonlinear
time
series and Volterra processes, an important class of nonlinear
processes.
\end{abstract}

%

\end{frontmatter}

\section{Introduction}\label{seintro}

Let $X_1, X_2, \ldots$ be independent, identically distributed
random variables with $\ep X_1=0$, $\ep X_1^2=1$. In their seminal
papers, Koml\'os, Major and Tusn\'ady (\citeyear{KomMajTus75,KomMajTus76}) proved that
under $\ep|X_1|^p<\infty$, $p>2$, there exists, after suitably
enlarging the probability space, a Wiener process $\{\B(t),
t\ge0\}$ such that, setting $S_n = \sum_{k=1}^n X_k$, we have
%
%
\begin{equation}
\label{kmt1} S_n=\B(n)+o\bigl(n^{1/p}\bigr) \qquad\mbox{a.s.}
\end{equation}
Assuming $\ep e^{t|X_1|}<\infty$ for some $t>0$, they obtained the
approximation
%
%
\begin{equation}
\label{kmt2} S_n=\B(n)+O(\log n) \qquad\mbox{a.s.}
\end{equation}
The remainder terms in (\ref{kmt1}) and (\ref{kmt2}) are optimal.
These results close a long development in probability theory
starting with the classical paper of Erd\H{o}s and Kac (\citeyear{ErdKac46})
introducing the method of \textit{invariance principle.} The ideas of
Erd\H{o}s and Kac were developed further by \citet{Doo49}, \citet{Don52},
\citet{Pro56} and others and led to the theory of weak
convergence of probability measures on metric spaces; see, for
example, \citet{Bil68}. In another direction, \citet{Str64} used the
Skorohod representation theorem to get an almost sure
approximation of partial sums of i.i.d. random variables by
Wiener process. Cs\"org\H{o} and R\'ev\'esz (\citeyear{CsoRev74}) showed that
using the quantile transform instead of Skorohod embedding yields
better approximation rates under higher moments and developing
this idea further,
Koml{\'o}s, Major and
Tusn{\'a}dy (\citeyear{KomMajTus75,KomMajTus76}) reached the final result in the
i.i.d. case.
Their results were extended to the independent, nonidentically
distributed case and for random variables taking values in
${\mathbb R}^d$, $d\ge2$, by Sakhanenko, Einmahl and Zaitsev; see
G\"otze and Zaitsev (\citeyear{GotZai08}) for history and references.

Due to the powerful consequences of KMT approximation [see, e.g., Cs\"
org\H{o} and Hall (\citeyear{CsoHal84}) or the books of Cs\"org\H{o} and
R\'ev\'esz (\citeyear{CsoRev81}) and \citet{ShoWel86} for the scope of
its applications], extending these results for dependent random
variables would have a great importance, but until recently,
little progress has been made in this direction. The dyadic
construction of Koml\'os, Major and Tusn\'ady is highly technical
and utilizes conditional large deviation techniques, which makes
it very difficult to extend to dependent processes. Recently a new
proof of the KMT result for the simple random walk via Stein's
method was given by \citet{Cha12}. The main motivation of his
paper was, as stated by the author, to get ``a more conceptual
understanding of the problem that may allow one to go beyond sums
of independent random variables.'' Using martingale approximation
and Skorohod embedding, \citet{ShaLu87} and \citet{Wu07} proved the
approximation
%
%
\begin{equation}
\label{wu2007} S_n=\sigma\B(n)+o\bigl(n^{1/p}(\log
n)^\gamma\bigr) \qquad\mbox{a.s.}
\end{equation}
with some $\sigma\ge0$, $\gamma>0$ for some classes of stationary
sequences $(X_k)$ satisfying $\ep X_1=0$, $\ep|X_1|^p<\infty$ for
some $2<p\le4$. \citet{LiuLin09} removed the logarithmic term
from (\ref{wu2007}), reaching the KMT bound $o(n^{1/p})$. Recently
\citet{MerRio12} and Dedecker, Doukhan and Merlev{\`e}de
(\citeyear{DedDouMer12}) extended these results for a much larger class of weakly
dependent processes. Note, however, that all
existing results in the dependent case concern the case $2\le p\le
4$ and the applied tools (e.g., Skorohod representation) limit the
accuracy of the approximation to $o(n^{1/4})$, regardless the
moment assumptions on $X_1$.

The purpose of the present paper is to develop a new approximation
technique enabling us to prove the KMT approximation (\ref{kmt1}) for
all $p>2$ and for a large class of dependent sequences.
Specifically, we will deal with stationary sequences allowing the
representation
%
%
\begin{equation}
\label{eqS41055} X_k=G(\ldots, \varepsilon_{k-1},
\varepsilon_k, \varepsilon_{k+1}, \ldots),\qquad k \in{\mathbb Z},
\end{equation}
where $\varepsilon_i$, $i \in\Z$, are i.i.d. random variables,
and $G\dvtx{\mathbb R}^{\mathbb Z}\to{\mathbb R}$ is a measurable
function. Sequences of this type have been studied intensively in
weak dependence theory [see, e.g., \citet{Bil68} or \citet{IbrLin71}], and many important time series models also have
a representation (\ref{eqS41055}). Processes of the type
(\ref{eqS41055}) also play an important role in ergodic theory,
as sequences generated by Bernoulli shift transformations. The
Bernoulli shift is a very important class of dynamical systems;
see \citet{Orn74} and \citet{Shi73} for the deep
Kolmogorov--Sinai--Ornstein isomorphism theory. There is a
substantial amount of research showing that various dynamical
systems are isomorphic to Bernoulli shifts. As a step further,
Weiss (\citeyear{Wei75}) asked,

\begin{quote}
``having shown that some physical system is Bernoullian, what does
that allow one to say about the system itself? To answer such
questions one must dig deeper and gain a better understanding of a
Bernoulli system.''
\end{quote}

\noindent Naturally, without additional assumptions one cannot
hope to prove KMT-type results (or even the CLT) for Bernoulli
systems; the representation (\ref{eqS41055}) allows stationary
processes that can exhibit a markedly non-i.i.d. behavior. For
limit theorems under dynamic assumptions, see \citet{HofKel82},
\citet{DenPhi84}, \citet{Den89}, Voln\'y (\citeyear{Vol99}),
\citet{MerRio12}. The classical approach to deal with
systems (\ref{eqS41055}) is to assume that $G$ is approximable
with finite dimensional functions in a certain technical sense;
see \citet{Bil68} or \citet{IbrLin71}. However,
this approach leads to a substantial loss of accuracy and does not
yield optimal results. In this paper we introduce a new, triadic
decomposition scheme enabling one to deduce directly, under the
dependence measure (\ref{eqpdm}) below, the asymptotic properties
of $X_n$ in (\ref{eqS41055}) from those of the~$\varepsilon_n$.
In particular, this allows us to carry over KMT approximation from
the partial sums of the $\varepsilon_n$ to those of $X_n$.

To state our weak dependence assumptions on the process in
(\ref{eqS41055}), assume $X_i \in{\cal L}^p$, $p
> 2$, namely $\| X_i \|_p := [\ep(|X_i|^p) ]^{1/p} < \infty$.
For $i\in\Z$ define the shift process ${\cal F}_i = (\varepsilon_{l+i},
l \in\Z)$. The central element of ${\cal F}_i$ (belonging to
$l=0$) is $\varepsilon_i$, and thus by (\ref{eqS41055}) we have
$X_i=G({\cal F}_i)$. Let $(\varepsilon'_j)_{j \in\Z}$ be an i.i.d.
copy of $(\varepsilon_j)_{j
\in\Z}$, and for $i, j \in\Z$ let ${\cal F}_{i, \{j\}}$ denote the process
obtained from ${\cal F}_i$ by replacing the coordinate $\varepsilon_j$ by
$\varepsilon_j'$. Put
%
%
\begin{equation}
\label{eqpdm} \delta_{i, p} = \|X_{i} - X_{i, \{0\}}
\|_p,\qquad \mbox{where } X_{i, \{0\}} = G({\cal F}_{i, \{0\}}).
\end{equation}
The above quantity can be interpreted as the dependence of $X_i$ on
$\varepsilon_0$ and $X_{i, \{0\}}$ is a coupled version of $X_{i}$
with $\varepsilon_0$ in the latter replaced by $\varepsilon'_0$. If
$G({\cal F}_i)$ does not functionally depend on $\varepsilon_0$,
then $\delta_{i, p} = 0$. Throughout the paper, for a random
variable $W = H({\cal F}_i)$, we use the notation $W_{ \{j\} } =
H({\cal F}_{i, \{j\} })$ for the $j$-coupled version of $W$.

The functional dependence measure (\ref{eqpdm}) is easy to work
with, and it is directly related to the underlying data-generating
mechanism. In our main result Theorem~\ref{thopsip}, we express our
dependence condition in terms of
%
%
\begin{equation}
\label{eqcpdm} \Theta_{i,p} = \sum_{|j|\ge i}
\delta_{j, p},\qquad i\ge0,
\end{equation}
which can be interpreted as the cumulative dependence of $(X_j)_{|j|
\ge i}$ on $\varepsilon_0$, or equivalently, the cumulative dependence
of $X_0$ on $\varepsilon_j$, $|j|\ge i$. Throughout the paper we
assume that the short-range dependence condition
%
%
\begin{equation}
\label{eqsrdpdm} \Theta_{0,p} < \infty
\end{equation}
holds. If (\ref{eqsrdpdm}) fails, then the process $(X_i)$ can be
long-range dependent, and the partial sum processes behave no longer
like Brownian motions. Our main result is introduced in Section~\ref{secmain}, where we also include some discussion on the
conditions. The proof is given in Section~\ref{secproof}, with the
proof of some useful lemmas postponed until Section~\ref{secuselem}.

\section{Main results}
\label{secmain} We introduce some notation. For $u \in\R$, let
$\lceil u \rceil=\break  \min\{ i \in\Z\dvtx i \ge u \}$ and $\lfloor u
\rfloor= \max\{ i \in\Z\dvtx i \le u \}$. Write the ${\cal L}^2$
norm $\| \cdot\| = \| \cdot\|_2$. Denote by ``$\Rightarrow$'' the
weak convergence. Before stating our main result, we first introduce
a central limit theorem for $S_n$. Assume that $X_i$ has mean zero,
$\ep(X_i^2) < \infty$, with covariance function $\gamma_i = \ep
(X_0 X_i)$, $i \in\Z$. Further assume that
%
%
\begin{equation}
\label{eqA181016} \sum_{i=-\infty}^\infty \bigl\|
\ep(X_i | {\cal G}_0) - \ep(X_i|{\cal
G}_{-1})\bigr\| < \infty,
\end{equation}
where ${\cal G}_i = (\ldots, \varepsilon_{i-1}, \varepsilon_i)$.
Then we have
%
%
\begin{equation}
\label{eqA181022} { {S_n} \over\sqrt n} \Rightarrow N\bigl(0,
\sigma^2\bigr) \qquad\mbox{where } \sigma^2 = \sum
_{i \in\Z} \gamma_i.
\end{equation}
Results of the above type have been known for several decades; see
\citet{Han79}, \citet{Woo92}, Voln\'y (\citeyear{Vol93}) and Dedecker and
Merlev\`ede (\citeyear{DedMer03}) among others. \citet{Wu05} pointed out the
inequality $\|\ep(X_i | {\cal G}_0) - \ep(X_i |{\cal G}_{-1})\| \le
\delta_{i,2}$. Hence~(\ref{eqA181016}) follows from $\Theta_{0, 2}
< \infty$. With stronger moment and dependence conditions, the
central limit theorem (\ref{eqA181022}) can be improved to strong
invariance principles.

There is a huge literature for central limit theorems and invariance
principles for stationary processes; see, for example, the monographs of
\citet{IbrLin71}, \citet{EbTa86}, \citet{Bra07},
\citet{Dedetal07} and \citet{Bil68}, among others. To establish
strong invariance principles, here we shall use the framework of stationary
process (\ref{eqS41055}) and its associated functional dependence
measures (\ref{eqpdm}). Many important processes in probability and statistics
assume this form; see the examples at the end of this section, where
also estimates for the
functional dependence measure $\delta_{i, p}$ are given.
The following theorem, which is the main result of our paper, provides optimal
KMT approximation for processes (\ref{eqS41055}) under suitable assumptions
on the functional dependence measure.

\begin{theorem}\label{thopsip}
Assume that $X_i \in{\cal L}^p$ with mean $0$, $p > 2$, and there
exists $\alpha> p$ such that
%
%
\begin{equation}
\label{eqsrdsip} \Xi_{\alpha, p}:= \sum_{j=-\infty}^\infty
|j|^{1/2 - 1/\alpha} \delta_{j, p}^{p/\alpha} < \infty.
\end{equation}
Further assume that there exists a positive integer sequence
$(m_k)_{k=1}^\infty$ such that
%
%
\begin{eqnarray}
\label{eqmk1} M_{\alpha, p}:= \sum_{k=1}^\infty3^{k - k\alpha/ p}
m_k^{\alpha/2
-1} &<& \infty,
\\
\label{eqmap} \sum_{k=1}^\infty
{ {3^{k p/2} \Theta_{m_k, p}^p} \over{3^k}} &<& \infty
\end{eqnarray}
and
%
%
\begin{equation}
\label{eqA18705p} \Theta_{m_k, p}+ \min_{l \ge0} \bigl(
\Theta_{l, p} + l 3^{k (2/p-1)}\bigr) = o \biggl(
{ {3^{k(1/p-1/2)}}\over{(\log k)^{1/2}}} \biggr).
\end{equation}
Then there exists a probability space $(\Omega_c, {\cal A}_c,
\pr_c)$ on which we can define random variables $X^{c}_i$ with the
partial sum process $S^{c}_{n}=\sum_{i=1}^{n} X^{c}_i$, and a
standard Brownian motion $\B_c(\cdot)$, such that
$(X^{c}_i)_{i \in{\mathbb Z}}
\eqd(X_i)_{i \in{\mathbb Z}}$ and
%
%
\begin{equation}
\label{eqsipA181103} S^{c}_{n} - \sigma
\B_c(n) = o_{ a.s.}\bigl(n^{1/p}\bigr) \qquad\mbox{in }
(\Omega_c, {\cal A}_c, \pr_c).
\end{equation}
\end{theorem}

Gaussian approximation results of type (\ref{eqsipA181103}) have
many applications in statistics. For example, \citet{WuZha07}
dealt with simultaneous inference of trends in time series. \citet{EubSpe93} considered a similar problem for independent
observations. As pointed out by and C. \citet{WuChiHoo98},
basic difficulties in the theory of simultaneous inference under
dependence are due to the lack of suitable Gaussian approximation.
Using a recent ``split'' form of approximation, Berkes, H\"ormann and
Schauer (\citeyear{BerHorSch11}) obtained asymptotic estimates for increments of
stationary processes with applications to change point tests.
Theorem \ref{thopsip} improves these results and provides
optimal rates. Many further
applications of the KMT theory for i.i.d. sequences also extend
easily for dependent samples via Theorem \ref{thopsip}.\vadjust{\goodbreak}

A crucial issue in applying Theorem \ref{thopsip} is to find the
sequence $m_k$ and to verify conditions (\ref{eqsrdsip}),
(\ref{eqmk1}), (\ref{eqmap}) and (\ref{eqA18705p}). If
$\Theta_{m, p}$ decays to zero at the rate $O(m^{-\tau} (\log
m)^{-A})$, where $\tau> 0$, then we have the following corollary.
An explicit form of $m_k$ can also be given. Let
%
%
\begin{equation}
\label{eqA19901p} \tau_p = { {p^2-4+(p-2)\sqrt{p^2 + 20p + 4}} \over{8 p}}.
\end{equation}

\begin{corollary}
\label{corsip} Assume that any one of the following holds:
\begin{longlist}[(iii)]
\item[(i)] $p > 4$ and $\Theta_{m, p} = O(m^{-\tau_p} (\log m)^{-A})$,
where $A > {2\over3} (1/p +1 + \tau_p)$;

\item[(ii)] $p = 4$ and $\Theta_{m, p} = O(m^{-1} (\log m)^{-A})$ with $A >
3/2$;

\item[(iii)] $2 < p < 4$ and $\Theta_{m, p} = O(m^{-1} (\log m)^{-1/p})$.
\end{longlist}

Then there exists $\alpha> p$ and an integer sequence
$m_k$ such that (\ref{eqsrdsip}), (\ref{eqmk1}), (\ref{eqmap})
and (\ref{eqA18705p}) are all satisfied. Hence the strong
invariance principle (\ref{eqsipA181103}) holds.
\end{corollary}

\begin{pf}If $\Theta_{m, p} = O(m^{-\tau} (\log m)^{-A})$, then
\begin{eqnarray*}
\Xi_{\alpha, p} &\le& \sum_{l=1}^\infty2^{l(1/2 - 1/\alpha)}
\sum_{j=2^{l-1}}^{2^l-1} \bigl(\delta_{j, p}^{p/\alpha}
+ \delta_{-j, p}^{p/\alpha}\bigr)
\\
&\le& \sum
_{l=1}^\infty2^{l(1/2 - 1/\alpha)} 2^{(l-1) (1-p/\alpha)} \Biggl(
\sum_{j=2^{l-1}}^{2^l-1} (\delta_{j, p} +
\delta_{-j, p}) \Biggr)^{p/\alpha}
\\
&\le& \sum
_{l=1}^\infty2^{l(3/2 - 1/\alpha-p/\alpha)} \Theta_{2^{l-1}, p}^{p/\alpha}
\\
&=& \sum_{l=1}^\infty2^{l(3/2 - 1/\alpha-p/\alpha)} O
\bigl[\bigl(2^{-l \tau} l^{-A}\bigr)^{p/\alpha}\bigr],
\end{eqnarray*}
which is finite if $3/2 < (1+p + p\tau) / \alpha$ or $3/2 = (1+p +
p\tau) / \alpha$ and $A p/\alpha> 1$.

(i) Write $\tau= \tau_p$. The quantity $\tau_p$ satisfies the
following equation:
%
%
\begin{equation}
\label{eqA17816} { {\tau-(1/2-1/p)} \over{\tau/p-1/4+1/(2p)}} = {2\over3}(1+p + p\tau).
\end{equation}
Let $\alpha= {2\over3} (1+p + p \tau_p)$. Then (\ref{eqsrdsip})
requires that $A p/\alpha> 1$, or $A > \alpha/ p$. Let
%
%
\begin{equation}
\label{eqA181015pm} m_k = \bigl\lfloor3^{ k(\alpha/p-1)/(\alpha/2-1)}
k^{-1/(\alpha/2-1)} (\log k)^{-1/(p/2-1)} \bigr\rfloor,
\end{equation}
which satisfies (\ref{eqmk1}). Then $\Theta_{m_k, p} =
O(m_k^{-\tau} k^{-A})$. If $A > \tau/ (\alpha/2 -1)$, then
(\ref{eqA18705p}) holds. If $A > \tau/ (\alpha/2 -1) + 1/ p$, then
(\ref{eqmap}) holds. Combining these three inequalities on $A$, we
have (i), since $\alpha/ p > \tau/ (\alpha/2 -1) + 1/ p$.

(ii) In this case we can choose $\alpha= 6$ and $m_k = \lfloor
3^{k/4} / k \rfloor$.

(iii) Since $2 < p < 4$, we can choose $\alpha$ such that
$(2+p)/(3-p/2) < \alpha< (2+4p)/3$ and $m_k = \lfloor3^{k
(1/2-1/p)} \log k \rfloor$.
\end{pf}

Corollary \ref{corsip} indicates that, to establish Gaussian
approximation for a Bernoulli shift process, one only needs to
compute the functional dependence measure $\delta_{i, p}$ in~(\ref{eqpdm}). In the following examples we shall deal with some
special Bernoulli process. Example~\ref{exnts} concerns some
widely used nonlinear time series, and Example~\ref{exVolterra}
deals with Volterra processes which play an important role in the
study of nonlinear systems.

\begin{example}
Consider the measure-preserving transformation $T x =  2 x\, \operatorname{mod}\, 1$ on $([0,1], {\cal B}, \pr)$, where $\pr$ is the Lebesgue
measure on $[0, 1]$. Let $U_0 \sim\operatorname{uniform} (0,1)$ have
the dyadic
expansion $U_0 = \sum_{j=0}^\infty\varepsilon_j / 2^{1+j}$, where
$\varepsilon_j$ are i.i.d. Bernoulli random variables with $\pr(
\varepsilon_j = 0) = \pr( \varepsilon_j = 1) = 1/2$. Then $U_i = T^i
U_0 = \sum_{j=i}^\infty\varepsilon_j / 2^{1+j-i}$, $i \ge0$; see
\citet{DenKel86} for a more detailed discussion. We now
compute the functional dependence measure for $X_i = g(U_i)$. Assume
that $\int_0^1 g(u) \,d u = 0$ and $\int_0^1 |g(u)|^p\, d u < \infty$,
$p > 2$. Then $\delta_{i,p} = 0$ if $i > 0$, and for $i \ge0$ we get
by stationarity
%
%
\begin{eqnarray}
\label{eqS5943} \delta_{-i,p}^p &=& \ep\bigl|g(U_0)
- g(U_{0, \{i\}})\bigr|^p
\nonumber
\\[-8pt]
\\[-8pt]
\nonumber
&=&{1\over2} \sum
_{j=1}^{2^i} \int_0^1
\biggl|g\biggl( { j \over{2^i}}+ {u \over{2^{i+1}}}\biggr) - g\biggl(
{ {j-1} \over{2^i}}+ {u \over{2^{i+1}}}\biggr)\biggr|^p\, d u.
\end{eqnarray}
If $X_i = g(U_i) = K(\sum_{j=i}^\infty a_{j-i} \varepsilon_j)$,
where $K$ is a Lipschitz continuous function and $\sum_{j=0}^\infty
|a_j| < \infty$, then $\delta_{i, p} = O(|a_i|)$. If $g$ has the Haar
wavelet expansion
%
%
\begin{equation}
\label{eqS71208} g(u) = \sum_{i=0}^\infty
\sum_{j=1}^{2^i} c_{i,j}
\phi_{i,j}(u),
\end{equation}
where $\phi_{i,j}(u) = 2^{i/2} \phi(2^i u-j)$ and $\phi(u) = \mathbf{
1}_{0\le u < 1/2} - \mathbf{ 1}_{1/2 \le u < 1}$, then for $i \ge0$,
%
%
\begin{equation}
\label{eqS71209} \delta_{-i, p}^p = O\bigl(2^{i(p/2-1)}
\bigr) \sum_{j=1}^{2^i} |c_{i,j}|^p.
\end{equation}
\end{example}

\begin{example}[(Nonlinear time series)]
\label{exnts} Consider the iterated
random function
%
%
\begin{equation}
\label{eqJ0210431} X_i = G(X_{i-1}, \varepsilon_i),
\end{equation}
where $\varepsilon_i$ are i.i.d. and $G$ is a measurable function
[\citet{DiaFre99}]. Many nonlinear time series including
ARCH, threshold autoregressive, random coefficient autoregressive and
bilinear autoregressive processes are of form (\ref{eqJ0210431}).
If there exists $p > 2$ and $x_0$ such that $G(x_0, \varepsilon_0)
\in{\cal L}^p$ and
%
%
\begin{equation}
\label{eqJ0210481} \ell_p = \sup_{x\not=x'}
{
{\| G(x, \varepsilon_0)-G(x', \varepsilon_0)\|_p}
\over{ |x-x'|}} < 1,
\end{equation}
then $\delta_{m, p} = O(\ell^m_p)$ and also $\Theta_{m, p} =
O(\ell^m_p)$ [\citet{WuSha04}]. Hence conditions in Corollary
\ref{corsip} are trivially satisfied, and thus (\ref{eqsipA181103})
holds.
\end{example}

\begin{example}
\label{exVolterra} In the study of nonlinear systems, Volterra
processes are of fundamental importance; see \citet{Sch80}, \citet{Rug81}, \citet{Cas85}, \citet{Pri88} and \citet{Ben90}, among
others. We consider the discrete-time process
%
%
\begin{equation}
\label{eqVolterra} X_n = \sum_{k=1}^\infty
\sum_{0 \le j_1 < \cdots< j_k} g_k(j_1, \ldots,
j_k) \varepsilon_{n-j_1} \cdots \varepsilon_{n-j_k},
\end{equation}
where $\varepsilon_i$ are i.i.d. with mean $0$, $\varepsilon_i \in
{\cal L}^p$, $p > 2$, and $g_k$ are called the $k$th order Volterra
kernel. Let
%
%
\begin{equation}
\label{eqQnk} Q_{n, k} = \sum_{n \in\{j_1, \ldots, j_k\},\ 0 \le j_1 < \cdots<
j_k}
g^2_k(j_1, \ldots, j_k).
\end{equation}
Assume for simplicity that $p$ is an even integer. Elementary
calculations show that there exists a constant $c_p$, only depending
on $p$, such that
%
%
\begin{equation}
\label{eqpdmVolterra} \delta_{n, p}^2 \le c_p
\sum_{k=1}^\infty\| \varepsilon_0
\|_p^{2k} Q_{n, k}.
\end{equation}
Assume that for some $\tau> 0$ and $A$,
%
%
\begin{equation}
\label{eqF51036}\qquad \sum_{k=1}^\infty\|
\varepsilon_0\|_p^{2k} \sum
_{j_k\ge m,\ 0 \le j_1 < \cdots< j_k} g^2_k(j_1, \ldots,
j_k) = O\bigl(m^{-1-2\tau}(\log m)^{-2A}\bigr)
\end{equation}
as $m \to\infty$. Then
%
%
\begin{equation}
\label{eqF51042} \sum_{n=m}^\infty
\delta_{n, p}^2 \le c_p \sum
_{k=1}^\infty\| \varepsilon_0
\|_p^{2k} \sum_{n=m}^\infty
Q_{n, k} = O\bigl(m^{-1-2\tau}(\log m)^{-2A}\bigr),
\end{equation}
which implies $\Theta_{m, p} = O(m^{-\tau} (\log m)^{-A})$ and
hence Corollary \ref{corsip} is applicable.
\end{example}

For further examples of processes allowing the representation (\ref{eqS41055}),
we refer to \citet{Wie58}, \citet{Ton90}, \citet{Pri88}, \citet{ShaWu07}, \citet{Wu11}
and the examples in Berkes, H\"ormann and Schauer (\citeyear{BerHorSch11}).

\section{\texorpdfstring{Proof of Theorem \protect\ref{thopsip}}{Proof of Theorem 2.1}}
\label{secproof} The proof of Theorem \ref{thopsip} is quite
intricate. To simplify the notation, we assume that $(X_i)$ is a
function of a one-sided Bernoulli shift,
%
%
\begin{equation}
\label{eqscp} X_{i} = G({\cal F}_i), \qquad\mbox{where } {\cal
F}_i = (\ldots, \varepsilon_{i-1}, \varepsilon_{i}),\vadjust{\goodbreak}
\end{equation}
where $\varepsilon_{k}, k \in\Z$, are i.i.d.
Clearly, in this case in (\ref{eqpdm})
we have $\delta_{i, p}=0$ for $i<0$.
As argued in \citet{Wu11}, (\ref{eqscp}) itself defines a very large
class of stationary
processes, and many widely used linear and nonlinear processes fall
within the framework of (\ref{eqscp}). Our argument can be extended
to the two-sided process (\ref{eqS41055}) in a straightforward
manner since our primary tool is the $m$-dependence approximation
technique. In Section~\ref{secTMB} we shall handle the
pre-processing work of truncation, $m$-dependence approximation and
blocking, and in Section~\ref{secGAA} we shall apply Sakhanenko's
(\citeyear{Sak06}) Gaussian approximation result to the transformed processes
and establish conditional Gaussian approximations. Section~\ref{secUGA} removes the conditioning, and an unconditional
Gaussian approximation is obtained. In Section~\ref{secA18RGA} we
refine the unconditional Gaussian approximation in Section~\ref{secUGA} by linearizing the variance function, so that one can
have the readily applicable form (\ref{eqsipA181103}).

\subsection{Truncation, $m$-dependence approximation and blocking}
\label{secTMB} For $a > 0$, define the truncation operator $T_a$ by
%
%
\begin{equation}
\label{eqtruncA13917} T_a(w) = \max\bigl(\min(w, a), -a\bigr), \qquad w
\in\R.
\end{equation}
Then $T_a$ is Lipschitz continuous and the Lipschitz constant is
$1$. For $n \ge2$ let $h_n = \lceil(\log n) / (\log3) \rceil$, so
that $3^{h_n-1} < n \le3^{h_n}$. Define
%
%
\begin{equation}
\label{eqA14725p} W_{k, l} = \sum_{i=1+3^{k-1}}^{l+3^{k-1}}
\bigl[T_{3^{k/p}} (X_i) - \ep T_{3^{k/p}}
(X_i)\bigr]
\end{equation}
and the $m_k$-dependent process
%
%
\begin{equation}
\label{eqtruncSA13933} \tilde X_{k,j} = \ep\bigl[ T_{3^{k/p}}
(X_j) | \varepsilon_{j-m_k}, \ldots, \varepsilon_{j-1},
\varepsilon_j\bigr] - \ep T_{3^{k/p}} (X_j).
\end{equation}
Let
%
%
\begin{equation}
\label{eqtruncSA13923} S_n^\dag= \sum
_{k=1}^{h_n-1} W_{k, 3^k-3^{k-1}} + \sum
_{i=1+3^{h_n-1}}^n \bigl[T_{3^{h_n/p}} (X_i)
- \ep T_{3^{h_n/p}} (X_i)\bigr]
\end{equation}
and
%
%
\begin{equation}
\label{eqA14814p} \tilde S_n = \sum_{k=1}^{h_n-1}
\tilde W_{k, 3^k-3^{k-1}} + \tilde W_{h_n, n - 3^{h_n-1}} \qquad\mbox{where } \tilde
W_{k, l} = \sum_{i=1+3^{k-1}}^{l+3^{k-1}} \tilde
X_{k,i}.
\end{equation}
If $n = 1$, we let $S_1^\dag= \tilde S_1 = 0$. Since $X_i \in{\cal
L}^p$, we have
%
%
\begin{equation}
\label{eqtruncSA13930} \max_{1\le i \le n} \bigl|S_i -
S_i^\dag\bigr| = o_\mathrm{ a.s.}\bigl(n^{1/p}\bigr).
\end{equation}
Note that there exists a constant $c_p$ such that, for all $k \ge
1$,
%
%
\begin{equation}
\label{eqA14820p} \Bigl\llVert \max_{3^{k-1} < l \le3^k} |\tilde
W_{k, l} - W_{k, l}| \Bigr\rrVert _p \le
c_p \bigl(3^k-3^{k-1}\bigr)^{1/2}
\Theta_{1+m_k, p}.
\end{equation}
Hence, by the Borel--Cantelli lemma and condition (\ref{eqmap}), we
have
%
%
\begin{equation}
\label{eqA14825p} \max_{1\le i \le n} \bigl|\tilde S_i -
S_i^\dag\bigr| = o_\mathrm{ a.s.}\bigl(n^{1/p}\bigr).
\end{equation}
Let $q_k = \lfloor2 \times3^{k-2} / m_k \rfloor- 2$. By
(\ref{eqmk1}), $m_k = o(3^{k(\alpha/p-1)/(\alpha/2-1)})$. Hence\break
$\lim_{k \to\infty} q_k = \infty$. Choose $K_0 \in\N$ such that
$q_k \ge2$ whenever $k \ge K_0$, and let $N_0 = 3^{K_0}$. For $k \ge
K_0$ define
%
%
\begin{equation}
\label{eqA14923p} B_{k, j} = \sum_{i = 1 + 3 j m_k + 3^{k-1}}^{ 3(j+1)m_k + 3^{k-1}}
\tilde X_{k,i},\qquad j=1, 2, \ldots, q_k.
\end{equation}
Let $B_{k, j} \equiv0$ if $k < K_0$. In the sequel we assume
throughout that $k \ge K_0$ and $n \ge N_0$. By Markov's inequality
and the stationarity of the process $(\tilde X_{k,i})_{i \in\Z}$,
%
%
\begin{eqnarray}
\label{eqA14927p}&& \pr \Biggl( \max_{1 \le l \le2 \times3^{k-1}} \Biggl\llvert \tilde
W_{k, l} - \sum_{j=1}^{\lfloor l/(3 m_k) \rfloor}
B_{k, j} \Biggr\rrvert \ge3^{k/p} \Biggr) \nonumber\\
&&\qquad\le
{ {2 \times3^{k-1}} \over{m_k}} \pr \Bigl( \max_{1 \le l \le3 m_k} | \tilde
W_{k, l}| \ge3^{k/p} \Bigr)
\\
&& \qquad\le{ {3^k \ep(\max_{1 \le l \le3 m_k}
| \tilde W_{k, l}|^\alpha)} \over{m_k 3^{k\alpha/p}}}.\nonumber
\end{eqnarray}
We define the functional dependence measure for the process
$(T_{3^{k/p}} (X_i))_{i \in\Z}$ as
%
%
\begin{equation}
\label{eqtruncSA14738} \delta_{k, j, \iota} = \bigl\|T_{3^{k/p}}
(X_i) - T_{3^{k/p}} (X_{i, \{i-j\}})\bigr\|_\iota,
\end{equation}
where $\iota\ge2$, and similarly the functional dependence measure
for $(\tilde X_{k, i})$ as
%
%
\begin{equation}
\tilde\delta_{k, j, \iota} = \|\tilde X_{k, i} - \tilde
X_{k, i, \{i-j\}}\|_\iota.
\end{equation}
For those dependence measures, we can easily have the following
simple relation:
%
%
\begin{equation}
\label{eqA181112} \tilde\delta_{k, j, \iota} \le\delta_{k, j, \iota},
\delta_{k, j, p} \le\delta_{j, p} \quad\mbox{and}\quad
\delta_{k, j, 2} \le\delta_{j, 2}.
\end{equation}
By the above relation, a careful check of the proof of Lemma
\ref{lemmomentA14} below indicates that,
under (\ref{eqsrdsip}) and
(\ref{eqmk1}), there exists a constant $c = c_{\alpha, p}$ such
that
%
%
\begin{equation}
\label{eqA14951p} \qquad\sum_{k=K_0}^\infty
{ {3^k}\over{m_k}} { {\ep(\max_{1\le l \le3 m_k} |\tilde W_{k, l}|^\alpha)}
\over{ 3^{k\alpha/p} }} \le c \bigl(M_{\alpha, p}
\Theta_{0, 2}^\alpha + \Xi_{\alpha, p}^\alpha +
\|X_1\|_p^p\bigr).
\end{equation}
The above inequality plays a critical role in our proof, and it will
be used again later. In (\ref{eqA14927p}), the largest index $j$ is
${\lfloor2 \times3^{k-1} /(3 m_k) \rfloor} = q_k + 2$. Note that
$B_{k, q_k}$ is independent of $B_{k+1, 1}$. This motivates us to
define the sum
%
%
\begin{equation}
\label{eqA141010p} \qquad S_n^\diamond= \sum
_{k=K_0}^{h_n-1} \sum_{j=1}^{q_k}
B_{k,j} + \sum_{j=1}^{\tau_n}
B_{h_n, j},\qquad \mbox{where } \tau_n = \biggl\lfloor
{ {n-3^{h_n-1}}
\over{3 m_{h_n}}} \biggr\rfloor- 2.\vadjust{\goodbreak}
\end{equation}
We emphasize that the sums $\sum_{j=1}^{q_k} B_{k,j}$, $k = 1, 2,
\ldots, h_n-1$ and $\sum_{j=1}^{\tau_n} B_{h_n, j}$ are mutually
independent. By (\ref{eqA14927p}), (\ref{eqA14951p}) and the
Borel--Cantelli lemma, we have
%
%
\begin{equation}
\label{eqA17109p} \max_{N_0 \le i \le n} \bigl|\tilde S_i -
S_i^\diamond\bigr| = o_\mathrm{ a.s.}\bigl(n^{1/p}\bigr),
\end{equation}
where we recall $N_0 = 3^{K_0}$. Summarizing the truncation
approximation (\ref{eqtruncSA13930}), the $m$-dependence
approximation (\ref{eqA14825p}) and the block approximation
(\ref{eqA17109p}), we have
%
%
\begin{equation}
\label{eqA20226} \max_{N_0 \le i \le n} \bigl|S_i -
S_i^\diamond\bigr| = o_\mathrm{ a.s.}\bigl(n^{1/p}\bigr),
\end{equation}
and by Lemma \ref{lemcnstctn} in Chapter~\ref{s4} it remains to show that
(\ref{eqsipA181103}) holds with~$S_n^\diamond$.

\subsection{Conditional Gaussian approximation}
\label{secGAA} For $3^{k-1} < i \le3^k$, $k \ge K_0$, let $G_k$ be
a measurable function such that
%
%
\begin{equation}
\label{eqmdrv} \tilde X_{k, i} = G_{k}(\varepsilon_{i - m_{k}},
\ldots, \varepsilon_i).
\end{equation}
Recall $q_k = \lfloor2 \times3^{k-2} / m_k \rfloor- 2$. For $j =
1, 2, \ldots, q_k$ define
%
%
\begin{equation}
\label{eqbbs1} {\cal J}_{k, j} = \bigl\{3^{k-1}+
(3j-1)m_k + l, l=1, 2, \ldots, m_k\bigr\}.
\end{equation}
Let $\mathbf{ a} = (\mathbf{ a}_{k, 3j}, 1\le j \le
q_k)_{k=K_0}^\infty
$ be
a vector of real numbers, where $\mathbf{ a}_{k, 3j} = (a_l, l \in
{\cal J}_{k, j})$, $j=1, \ldots, q_k$. Define the random functions
\begin{eqnarray*}
F_{k, 3j}(\mathbf{a}_{k, 3j}) &=& \sum
_{i=1+(3j-1)m_k}^{3j m_k} G_k(a_{i+3^{k-1}},
\ldots,a_{3 j m_k+3^{k-1}},
\\
& & \hspace*{77pt}\varepsilon_{3j m_k+1+3^{k-1}},\ldots,
\varepsilon_{i+m_k+3^{k-1}});
\\
F_{k,3j+1}&=&\sum
_{i=1+3jm_k}^{(3j+1)m_k} G_k(\varepsilon_{i+3^{k-1}},
\ldots, \varepsilon_{(3j+1)m_k+3^{k-1}},
\\
& &\hspace*{60pt} \varepsilon_{(3j+1)m_k+1+3^{k-1}},\ldots,
\varepsilon_{i+m_k+3^{k-1}});
\\
F_{k,3j+2}(\mathbf{a}_{k,3j+3})
&=&\sum_{i=1+(3j+1)m_k}^{(3j+2)m_k} G_k(
\varepsilon_{i+3^{k-1}},\ldots,\varepsilon_{(3j+2)m_k+3^{k-1}},
\\
&&\hspace*{78pt}
a_{(3j+2)m_k+1+3^{k-1}},\ldots,a_{i+m_k+3^{k-1}}).
\end{eqnarray*}
Let $\etav_{k, 3j} = (\varepsilon_l, l \in{\cal J}_{k, j})$,
$j=1, \ldots, q_k$, and $\etav= (\etav_{k, 3j}, 1\le j \le
q_k)_{k=K_0}^\infty$. Then
%
%
\begin{equation}
\label{eqA16744} B_{k, j} = F_{k, 3j}(\etav_{k, 3j}) +
F_{k,3j+1} + F_{k,3j+2}(\etav_{k, 3j+3}).
\end{equation}
Note that $\ep F_{k,3j+1} = 0$. Define the mean functions
\[
\Lambda_{k, 0} (\mathbf{a}_{k, 3j}) = \ep F_{k, 3j}(
\mathbf{a}_{k,
3j}),\qquad \Lambda_{k, 2} (\mathbf{a}_{k, 3j+3})
= \ep F_{k,3j+2} (\mathbf{a}_{k,3j+3}).
\]
Introduce the centered process
%
%
\begin{eqnarray}
\label{eqA161219} \qquad Y_{k, j}(\mathbf{a}_{k, 3j},
\mathbf{a}_{k, 3j+3}) &=& \bigl[F_{k, 3j}(\mathbf{a}_{k, 3j})
- \Lambda_{k,0}(\mathbf{a}_{k, 3j})\bigr]
\nonumber
\\[-8pt]
\\[-8pt]
\nonumber
&& {}+
F_{k, 3j+1} + \bigl[F_{k, 3j+2}(\mathbf{a}_{k, 3j+3}) -
\Lambda_{k,2}(\mathbf{a}_{k, 3j+3})\bigr].
\end{eqnarray}
Then $Y_{k, j}(\mathbf{a}_{k, 3j}, \mathbf{a}_{k, 3j+3})$, $j = 1,
\ldots, q_k$, $k \ge K_0$, are mean zero independent random
variables with variance function
%
%
\begin{eqnarray}
\label{eqA16752} V_k(\mathbf{a}_{k, 3j},
\mathbf{a}_{k, 3j+3}) &=& \bigl\|Y_{k, j}(\mathbf{a}_{k, 3j},
\mathbf{a}_{k, 3j+3})\bigr\|^2
\nonumber\\
&=& \bigl\|F_{k, 3j}(
\mathbf{a}_{k, 3j}) - \Lambda_{k,0}(\mathbf{a}_{k, 3j})
\bigr\|^2 + \|F_{k, 3j+1}\|^2
\nonumber\\
&&{} + 2 \ep\bigl
\{F_{k, 3j+1} \bigl[F_{k, 3j}(\mathbf{a}_{k, 3j}) -
\Lambda_{k,0}(\mathbf{a}_{k, 3j})\bigr] \bigr\}
\\
&&{} + \bigl\|
F_{k, 3j+2}(\mathbf{a}_{k, 3j+3}) - \Lambda_{k,2}(
\mathbf{a}_{k, 3j+3}) \bigr\|^2\nonumber
\\
&&{} + 2 \ep\bigl
\{F_{k, 3j+1} \bigl[F_{k, 3j+2}(\mathbf{a}_{k, 3j+3}) -
\Lambda_{k,2}(\mathbf{a}_{k, 3j+3})\bigr] \bigr\},\nonumber
\end{eqnarray}
since $[F_{k, 3j}(\mathbf{a}_{k, 3j}) - \Lambda_{k,0}(\mathbf{a}_{k,
3j})]$ and $[F_{k, 3j+2}(\mathbf{a}_{k, 3j+3}) -
\Lambda_{k,2}(\mathbf{a}_{k, 3j+3})]$ are independent. Following the
definition of $S_n^\diamond$ in (\ref{eqA141010p}), we let
%
%
\begin{eqnarray}
\label{eqF281} H_n(\mathbf{ a})& =& \sum_{k=K_0}^{h_n-1}
\sum_{j=1}^{q_k} Y_{k, j}(
\mathbf{a}_{k, 3j}, \mathbf{a}_{k, 3j+3})
\nonumber
\\[-8pt]
\\[-8pt]
\nonumber
&&{} + \sum
_{j=1}^{\tau_n} Y_{h_n, j}(\mathbf{a}_{h_n, 3j},
\mathbf{a}_{h_n, 3j+3}).
\end{eqnarray}
Define the mean function
\begin{eqnarray*}
M_n(\mathbf{ a}) &=& \sum_{k=K_0}^{h_n-1}
\sum_{j=1}^{q_k} \bigl[\Lambda_{k, 0}
(\mathbf{a}_{k, 3j}) + \Lambda_{k, 2} (\mathbf{a}_{k, 3j+3})
\bigr]
\\
& & {}+ \sum_{j=1}^{\tau_n} \bigl[
\Lambda_{h_n, 0} (\mathbf{a}_{h_n, 3j}) + \Lambda_{h_n, 2} (
\mathbf{a}_{h_n, 3j+3})\bigr],
\end{eqnarray*}
and the variance of $H_n(\mathbf{ a})$,
\[
Q_n(\mathbf{ a}) = \sum_{k=K_0}^{h_n-1}
\sum_{j=1}^{q_k} V_{k} (
\mathbf{a}_{k, 3j}, \mathbf{a}_{k, 3j+3}) + \sum
_{j=1}^{\tau_n} V_{h_n} (\mathbf{a}_{h_n, 3j},
\mathbf{a}_{h_n, 3j+3}).
\]
Let
%
%
\begin{eqnarray}
\label{eqA16818} V^\circ_k(\mathbf{a}_{k, 3j}) &=&
\bigl\|\bigl[F_{k, 3j}(\mathbf{a}_{k, 3j}) - \Lambda_{k,0}(
\mathbf{a}_{k, 3j})\bigr]
\nonumber\\
&&\hspace*{6pt}{} + F_{k, 3j+1} +
\bigl[F_{k, 3j+2}(\mathbf{a}_{k, 3j}) - \Lambda_{k,2}(
\mathbf{a}_{k, 3j})\bigr]\bigr\|^2
\nonumber\\
&=& \bigl\|F_{k, 3j}(
\mathbf{a}_{k, 3j}) - \Lambda_{k,0}(\mathbf{a}_{k, 3j})
\bigr\|^2 + \| F_{k, 3j+1}\|^2
\nonumber\\
&&{} + 2\ep\bigl
\{F_{k, 3j+1} \bigl[F_{k, 3j}(\mathbf{a}_{k, 3j}) -
\Lambda_{k,0}(\mathbf{a}_{k, 3j})\bigr]\bigr\}
\nonumber\\
&&{} +
\bigl\|F_{k, 3j+2}(\mathbf{a}_{k, 3j}) - \Lambda_{k,2}(
\mathbf{a}_{k, 3j})\bigr\|^2
\\
&&{} + 2\ep\bigl\{F_{k, 3j+1}
\bigl[F_{k, 3j+2}(\mathbf{a}_{k, 3j}) - \Lambda_{k,2}(
\mathbf{a}_{k, 3j})\bigr]\bigr\},
\nonumber\\
L_k(
\mathbf{a}_{k, 3j}) &=& \bigl\|F_{k, 3j+1} + \bigl[F_{k, 3j+2}(
\mathbf{a}_{k, 3j}) - \Lambda_{k,2}(\mathbf{a}_{k, 3j})
\bigr]\bigr\|^2
\nonumber\\
&=& \bigl\|F_{k, 3j+1}\bigr\|^2 + \bigl\|
\bigl[F_{k, 3j+2}(\mathbf{a}_{k, 3j}) - \Lambda_{k,2}(
\mathbf{a}_{k, 3j})\bigr]\bigr\|^2
\nonumber\\
&&{} + 2 \ep\bigl\{
F_{k, 3j+1} \bigl[F_{k, 3j+2}(\mathbf{a}_{k, 3j}) -
\Lambda_{k,2}(\mathbf{a}_{k, 3j})\bigr] \bigr\}.\nonumber
\end{eqnarray}
By the formulas of $V_k(\mathbf{a}_{k, 3j}, \mathbf{a}_{k, 3j+3})$
in (\ref{eqA16752}) and $V^\circ_k(\mathbf{a}_{k, 3j})$ and
$L_k(\mathbf{a}_{k, 3j})$ in~(\ref{eqA16818}), we have the
following identity:
%
%
\begin{equation}
\label{eqA16819} L_k(\mathbf{a}_{k, 3}) + \sum
_{j=1}^t V_k(\mathbf{a}_{k, 3j},
\mathbf{a}_{k, 3j+3}) = \sum_{j=1}^t
V^\circ_k(\mathbf{a}_{k, 3j}) + L_k(
\mathbf{a}_{k, 3+3t})
\end{equation}
holds for all $t \ge1$. The above identity motivates us to
introduce the auxiliary process
%
%
\begin{equation}
\label{eqA16835} \Gamma_n(\mathbf{a}) = \sum
_{k=K_0}^{h_n-1} L_k(\mathbf{a}_{k, 3})^{1/2}
\zeta_k + L_{h_n}(\mathbf{a}_{h_n, 3})^{1/2}
\zeta_{h_n},
\end{equation}
where $\zeta_l, l \in\Z$, are i.i.d. standard normal random
variables which are independent of $(\varepsilon_i)_{i \in\Z}$.
Then in view of (\ref{eqA16819}), the variance of $H_n(\mathbf{ a}) +
\Gamma_n(\mathbf{a})$ is given by
%
%
\begin{eqnarray}
\label{eqA16840} Q^\circ_n(\mathbf{ a}) &=& \sum
_{k=K_0}^{h_n-1} \Biggl[\sum_{j=1}^{q_k}
V^\circ_k(\mathbf{a}_{k, 3j}) + L_k(
\mathbf{a}_{k, 3+3 q_k}) \Biggr]
\nonumber
\\[-8pt]
\\[-8pt]
\nonumber
&& {}+ \sum_{j=1}^{\tau_n}
\bigl[ V^\circ_{h_n} (\mathbf{a}_{h_n, 3j}) +
L_{h_n}(\mathbf{a}_{h_n, 3+3 \tau_n}) \bigr].
\end{eqnarray}
In studying $H_n(\mathbf{ a}) + \Gamma_n(\mathbf{a})$, for notational
convenience, for $j = 0$ we let $Y_{k,0}(\mathbf{ a}_{k, 0},  \mathbf{
a}_{k, 3}) = L_k(\mathbf{a}_{k, 3})^{1/2} \zeta_k$. We shall now
apply Sakhanenko's (\citeyear{Sak91,Sak06}) Gaussian approximation result. To
this end, for $x > 0$, we define
%
%
\begin{eqnarray}
\label{eqA161125}&& \Psi_h(\mathbf{ a}, x, \alpha) \nonumber\\
&&\qquad= \sum
_{k=K_0}^h \sum_{j=0}^{q_k}
\ep\min\bigl\{ \bigl|Y_{k,j}(\mathbf{ a}_{k, 3j}, \mathbf{
a}_{k, 3j+3}) / x\bigr|^\alpha, \bigl|Y_{k,j}(\mathbf{
a}_{k, 3j}, \mathbf{ a}_{k, 3j+3}) / x\bigr|^2 \bigr\}
\\
&&\qquad
\le \sum_{k=K_0}^h \sum
_{j=0}^{q_k} \ep\bigl|Y_{k,j}(\mathbf{
a}_{k, 3j}, \mathbf{ a}_{k, 3j+3}) / x\bigr|^\alpha.\nonumber
\end{eqnarray}
By Theorem 1 in \citet{Sak06}, there exists a probability space
$(\Omega_\mathbf{ a}, {\cal A}_\mathbf{ a}, \pr_\mathbf{ a})$ on which
we can
define a standard Brownian motion $\B_\mathbf{ a}$ and random variables
$R_{k, j}^\mathbf{ a}$ such that the distributional equality
%
%
\begin{equation}
\label{eqeqRy} \bigl(R_{k,j}^\mathbf{ a}\bigr)_{0 \le j \le q_k, k \ge K_0}
\stackrel{\cal D}{=} \bigl(Y_{k,j}(\mathbf{ a}_{k, 3j}, \mathbf{
a}_{k, 3j+3})\bigr)_{0 \le j \le q_k, k \ge K_0}
\end{equation}
holds, and, for the partial sum processes
%
%
\begin{equation}
\label{eqA161218} \qquad\Upsilon_n^\mathbf{ a} = \sum
_{k=K_0}^{h-1} \sum_{j=1}^{q_k}
R_{k,j}^\mathbf{ a} + \sum_{j=1}^{\tau_n}
R_{h_n,j}^\mathbf{ a} \qquad\mbox{and}\qquad \mu_n^\mathbf{
a} = \sum_{k=K_0}^{h-1} R_{k,0}^\mathbf{
a} + R_{h_n,0}^\mathbf{ a},
\end{equation}
we have for all $x > 0$ and $\alpha> p$ that
%
%
\begin{equation}
\label{eqscha7} \pr_\mathbf{ a} \Bigl[ \max_{N_0 \le i \le3^h}
\bigl\llvert \bigl(\Upsilon_i^\mathbf{ a}+
\mu_i^\mathbf{ a}\bigr) - \B_\mathbf{ a}\bigl(
Q^\circ_i(\mathbf{ a}) \bigr) \bigr\rrvert \ge
c_0 \alpha x \Bigr] \le\Psi_h(\mathbf{ a}, x, \alpha).
\end{equation}
Here $c_0$ is an absolute constant. By Jensen's inequality, for both
$j=0$ and $j > 0$, there exists a constant $c_\alpha$ such that
%
%
\begin{equation}
\ep\bigl[\bigl|Y_{k,j}(\etav_{k, 3j}, \etav_{k, 3j+3})\bigr|^\alpha
\bigr] \le c_\alpha\ep\bigl(|\tilde W_{k, m_k}|^\alpha
\bigr).
\end{equation}
In (\ref{eqscha7}) we let $x = 3^{h/p}$ and by Lemma \ref{lemasct}
in the next chapter [see also (\ref{eqA14951p})],
%
%
\begin{eqnarray}
\label{eqA15830} \sum_{h=K_0}^\infty\ep\bigl[
\Psi_h\bigl(\etav, 3^{h/p}, \alpha\bigr)\bigr] &\le& \sum
_{h=K_0}^\infty\sum
_{k=K_0}^h { {q_k+1} \over{3^{\alpha h/p}}}
c_\alpha\ep\bigl(|\tilde W_{k, m_k}|^\alpha\bigr)
\nonumber\\
&
\le& \sum_{k=K_0}^\infty\sum
_{h=k}^\infty { {3^k c_\alpha} \over{m_k 3^{\alpha h/p}}} \ep \Bigl(\max
_{1\le l \le3 m_k} |\tilde W_{k, l}|^\alpha \Bigr)
\\
&
<& \infty.\nonumber
\end{eqnarray}
Hence, by the Borel--Cantelli lemma, we obtain
%
%
\begin{equation}
\label{eqA161140} \max_{i \le n}\bigl |\bigl(\Upsilon_i^{\etav}+
\mu_i^{\etav}\bigr) - \B_{\etav}\bigl(Q^\circ_i({
\etav})\bigr)\bigr| = o_\mathrm{ a.s.}\bigl(n^{1/p}\bigr).
\end{equation}
The probability space for the above almost sure convergence is
%
%
\begin{equation}
\label{eqps3} (\Omega_*, {\cal A}_*, \pr_*) = (\Omega, {\cal A}, \pr) \times\prod
_{\tau\in\Omega} (\Omega_{\etav(\tau)}, {\cal
A}_{\etav(\tau)}, \pr_{\etav(\tau)}),
\end{equation}
where $(\Omega, {\cal A}, \pr)$ is the probability space on which
the random variables $(\varepsilon_i)_{i \in\Z}$ are defined and,
for a set $A \subset\Omega_*$ with $A \in{\cal A}_*$, the
probability measure $\pr_*$ is defined as
%
%
\begin{equation}
\label{eq} \pr_*(A) = \int_\Omega\pr_{\etav(\omega)}
(A_\omega) \pr(d \omega),
\end{equation}
where $A_\omega$ is the $\omega$-section of $A$. Here we recall
that, for each $\mathbf{ a}$, $(\Omega_\mathbf{ a}, {\cal A}_\mathbf
{ a},
\pr_\mathbf{ a})$ is the probability space carrying $\B_\mathbf{
a}$ and
$R_{k, j}^\mathbf{ a}$ given $\etav= \mathbf{ a}$. On the
probability space
$(\Omega_*, {\cal A}_*, \pr_*)$, the random variable
$R_{k,j}^{\etav}$ is defined as $R_{k, j}^{\etav} (\omega,
\theta(\cdot)) = R_{k, j}^{\etav(\omega)}(\theta(\omega))$, where
$(\omega,\theta(\cdot)) \in\Omega_*$, $\theta(\cdot)$ is an element
in $\prod_{\tau\in\Omega} \Omega_{\etav(\tau)}$ and $\theta
(\tau)
\in\Omega_{\etav(\tau)}$, $\tau\in\Omega$. The other random
processes $\mu_i^{\etav}$ and $\B_{\etav}(Q^\circ_i({\etav}))$ can
be similarly defined.

\subsection{Unconditional Gaussian approximation}
\label{secUGA} In this subsection we shall work with the processes
$\Upsilon_i^{\etav}$, $\mu_i^{\etav}$ and
$\B_{\etav}(Q^\circ_i({\etav}))$.
Based on (\ref{eqA16840}), we
can construct i.i.d. standard normal random variables $Z^{\mathbf{ a}}_{i, l}, i, l \in \Z$, and standard normal random variables
$\mathcal{ G}^{\mathbf{ a}}_{i, l}$, such that
%
%
\begin{equation}
\label{eqA161157} \B_\mathbf{ a}\bigl( Q^\circ_n(
\mathbf{ a})\bigr) = \varpi_n(\mathbf{a}) + \varphi_n(
\mathbf{a}),
\end{equation}
where
\begin{eqnarray*}
\varpi_n(\mathbf{a}) &=& \sum_{k=K_0}^{h_n-1}
\sum_{j=1}^{q_k} V^\circ_k(
\mathbf{a}_{k, 3j})^{1/2} Z^\mathbf{ a}_{k, j} +
\sum_{j=1}^{\tau_n} V^\circ_{h_n}
(\mathbf{a}_{h_n, 3j})^{1/2} Z^\mathbf{
a}_{h_n, j},
\\
\varphi_n(\mathbf{a}) &=& \sum
_{k=K_0}^{h_n-1} L_k(\mathbf{a}_{k, 3+3 q_k})^{1/2}
{\mathcal G}^{\mathbf{a}}_{k, 1+q_k} + L_{h_n}(\mathbf{a}_{h_n, 3+3 \tau_n})^{1/2}
{\mathcal G}^\mathbf{ a}_{h_n, 1+\tau_n}.
\end{eqnarray*}
%
In particular,
\begin{eqnarray*}
{V^\circ_{h_n} (\mathbf{a}_{h_n, 3j})^{1/2}}
Z^{\mathbf {a}}_{h_n, j} &=& \B_{\mathbf {a}}\Biggl( Q^\circ_{3^{h_n-1}}(\mathbf {a}) +
\sum_{j'=1}^j V^\circ_{h_n} (\mathbf{a}_{h_n, 3j'})\Biggr) \\
&&{}- \B_{\mathbf{ a}}\Biggl(
Q^\circ_{3^{h_n-1}}({\mathbf {a}}) + \sum_{j'=1}^{j-1} V^\circ_{h_n}
(\mathbf{a}_{h_n, 3j'})\Biggr)
\end{eqnarray*}
and
\[
L_{h_n}(\mathbf{a}_{h_n, 3+3
\tau_n})^{1/2} {\cal G}^{\bf a}_{h_n, 1+\tau_n} = \B_{\mathbf a}\bigl(
Q^\circ_n({\mathbf a})\bigr) -  \B_{\mathbf a}\Biggl( Q^\circ_{3^{h_n-1}}({\bf a}) +
\sum_{j=1}^{\tau_n} V^\circ_{h_n} (\mathbf{a}_{h_n, 3j})\Biggr).
\]
Note that the standard normal random variables ${\cal G}^{\mathbf a}_{i,
l}, i,l,$ can be possibly dependent and $({\cal G}^{\bf a}_{i,
l})_{i l}$ and $(Z^{\mathbf a}_{i, l})_{i l}$ can also be possibly
dependent.

Let $Z^\star_{i, l}, i, l \in\Z$, independent of
$(\varepsilon_j)_{j \in\Z}$, be also i.i.d. standard normal random
variables, and define
\[
\Phi_n =\sum_{k=K_0}^{h_n-1}
\sum_{j=1}^{q_k} V^\circ_k(
\etav_{k, 3j})^{1/2} Z^\star_{k, j} + \sum
_{j=1}^{\tau_n} V^\circ_{h_n}
(\etav_{h_n, 3j})^{1/2} Z^\star_{h_n, j}.
\]
%
Since $Z^\mathbf{ a}_{i, l}$, are i.i.d. standard normal, the conditional
distribution $[\varpi_n(\etav) | \etav= \mathbf{ a}]$, namely the
distribution of $\varpi_n(\mathbf{a})$, is same as that of $\Phi_n$.
Hence
%
%
\begin{equation}
(\Phi_i)_{i \ge N_0}
 \stackrel{\cal D}{=} \bigl(\varpi_i(\etav)\bigr)_{i \ge N_0}.
\end{equation}
By Jensen's inequality, $\ep[|L_k(\etav_{k, 3j+3})^{1/2}|^\alpha]
\le3^{\alpha} \ep(|\tilde W_{k, m_k}|^\alpha)$. By
(\ref{eqA14951p}),
%
%
\begin{eqnarray}
\label{eqA17228}&& \sum_{k=K_0}^\infty\pr \Bigl(
\max_{1\le j \le q_k} \bigl|L_k(\etav_{k, 3j+3})^{1/2}
{\mathcal G}_{k, 1+j}^{\etav} \bigr| \ge3^{k/p} \Bigr)\nonumber\\
&&\qquad \le \sum
_{k=K_0}^\infty q_k
{ {\ep[|L_k(\etav_{k, 3})^{1/2} {\mathcal G}_{k, 1}^{\etav} |^\alpha]}
\over{3^{k\alpha/p}} }
\nonumber
\\[-8pt]
\\[-8pt]
\nonumber
&&\qquad\le \sum_{k=K_0}^\infty
q_k { {c_\alpha\ep(|\tilde W_{k, m_k}|^\alpha)}
\over{3^{k\alpha/p}} }
\\
&&\qquad <  \infty,\nonumber
\end{eqnarray}
which by the Borel--Cantelli lemma implies
%
%
\begin{equation}
\label{eqA17203} \max_{i\le n}\bigl |\varphi_i(\etav)\bigr| =
o_\mathrm{ a.s.}\bigl(n^{1/p}\bigr).
\end{equation}
The same argument also implies that $\max_{i\le n} |\Gamma_i(\etav)|
= o_\mathrm{ a.s.}(n^{1/p})$ and consequently
%
%
\begin{equation}
\label{eqA161200} \max_{i\le n}\bigl |\mu_i^{\etav}\bigr|
= o_\mathrm{ a.s.}\bigl(n^{1/p}\bigr)
\end{equation}
in view of (\ref{eqeqRy}) with $j = 0$. Hence by
(\ref{eqA161140}) and (\ref{eqA161157}), we have $ \max_{i \le n}
|\Upsilon_i^{\etav}-\varpi_i(\etav)| = o_\mathrm{ a.s.}(n^{1/p})$.
Observe that, by (\ref{eqeqRy}), (\ref{eqA161218}),
(\ref{eqA16744}) and (\ref{eqA161219}), we have the distributional
equality
%
%
\begin{equation}
\label{eqA161220} \bigl(\Upsilon_i^{\etav} + M_i(
\etav)\bigr)_{i \ge N_0} \stackrel{\cal D}{=} \bigl(S_i^\diamond
\bigr)_{i \ge N_0},
\end{equation}
where we recall (\ref{eqA141010p}) for the definition of
$S_n^\diamond$. Then it remains to establish a strong invariance
principle for $\Phi_n + M_n(\etav)$. To this end, let
%
%
\begin{equation}
A_{k, j} = V^\circ_k(\etav_{k, 3j})^{1/2}
Z^\star_{k, j} + \Lambda_{k, 0} (
\etav_{k, 3j}) + \Lambda_{k, 2} (\etav_{k, 3j}),
\end{equation}
which are independent random variables for $j = 1, \ldots, q_k$ and
$k \ge K_0$, and let
%
%
\begin{equation}
\label{eqA17220p} S_n^\natural= \sum
_{k=K_0}^{h_n-1} \sum_{j=1}^{q_k}
A_{k,j} + \sum_{j=1}^{\tau_n}
A_{h_n, j}
\end{equation}
and $R_n^\natural= \Phi_n + M_n(\etav) - S_n^\natural$. Note that
\[
R_n^\natural= \sum_{k=K_0}^{h_n-1}
\bigl[\Lambda_{k, 2} (\etav_{k, 3+3q_k}) - \Lambda_{k, 2} (
\etav_{k, 3})\bigr] + \bigl[\Lambda_{h_n, 2} (
\etav_{h_n, 3+3 \tau_n}) - \Lambda_{h_n, 2} (\etav_{h_n, 3})\bigr].
\]
Then using the same argument as in (\ref{eqA17228}), we have
%
%
\begin{equation}
\label{eqA17231p} \max_{i \le n}\bigl |R_i^\natural\bigr|
= \max_{i \le n} \bigl|\Phi_i + M_i(\etav) -
S_i^\natural\bigr| = o_\mathrm{ a.s.}\bigl(n^{1/p}\bigr).
\end{equation}
The variance of $S_n^\natural$ equals to
%
%
\begin{eqnarray}
\label{eqA20212p} \sigma_n^2 &=& \sum
_{k=K_0}^{h_n-1} \sum_{j=1}^{q_k}
\|A_{k, j}\|^2 + \sum_{j=1}^{\tau_n}
\|A_{h_n, j}\|^2
\nonumber
\\[-8pt]
\\[-8pt]
\nonumber
& =& \sum_{k=K_0}^{h_n-1}
q_k \|A_{k, 1}\|^2 + \tau_n
\|A_{h_n, 1}\|^2.
\end{eqnarray}
Again by Theorem 1 in \citet{Sak06}, on the same probability
space that defines $(A_{k, j})_{1\le j \le q_k, k \ge K_0}$, by
the argument in (\ref{eqscha7})--(\ref{eqA161140}), there exists a
standard Brownian motion $\B$ such that
%
%
\begin{equation}
\label{eqA161249} \max_{i\le n} \bigl|S_i^\natural-
\B\bigl(\sigma_i^2\bigr)\bigr| = o_\mathrm{ a.s.}
\bigl(n^{1/p}\bigr).
\end{equation}

\subsection{Regularizing the Gaussian approximation}
\label{secA18RGA} In this section we shall regularize the Gaussian
approximation (\ref{eqA161249}) by replacing the variance function
$\sigma_i^2$ by the asymptotic linear form $\phi_i$ or the linear
form $i \sigma^2$, and the latter is more easily usable. By
(\ref{eqA16818}), we obtain
%
%
\begin{eqnarray}
V^\circ_k(\mathbf{a}_{k, 3j}) &=&
\bigl\|F_{k, 3j}(\mathbf{a}_{k, 3j})\bigr\|^2 -
\Lambda_{k,0}(\mathbf{a}_{k, 3j})^2 + \|
F_{k, 3j+1}\|^2
\nonumber\\
&&{} + 2\ep\bigl\{F_{k, 3j+1}
F_{k, 3j}(\mathbf{a}_{k, 3j}) \bigr\}
\nonumber
\\[-8pt]
\\[-8pt]
\nonumber
&&{} +
\bigl\|F_{k, 3j+2}(\mathbf{a}_{k, 3j})\bigr\|^2 -
\Lambda_{k,2}(\mathbf{a}_{k, 3j})^2
\\
&&{} + 2\ep
\bigl\{F_{k, 3j+1} F_{k, 3j+2}(\mathbf{a}_{k, 3j}) \bigr\},\nonumber
\end{eqnarray}
which, by the expression of $A_{k, j}$, implies that
%
%
\begin{eqnarray}
\|A_{k, j}\|^2 &=& \ep\bigl[V^\circ_k(
\etav_{k, 3j})\bigr] + \ep\bigl[\Lambda_{k, 0} (
\etav_{k, 3j}) + \Lambda_{k, 2} (\etav_{k, 3j})
\bigr]^2
\nonumber
\\[-8pt]
\\[-8pt]
\nonumber
& =& 3 \ep\bigl[\tilde W_{k, m_k}^2 +
2 \tilde W_{k, m_k} (\tilde W_{k, 2 m_k}-\tilde W_{k, m_k})
\bigr].
\end{eqnarray}
Let $\tilde\gamma_{k, i} = \ep(\tilde X_{k,0} \tilde X_{k,i} )$.
Then $\nu_k := \|A_{k, j}\|^2 / (3 m_k)$ has the expression
%
%
\begin{eqnarray}
\label{eqA19513p} \nu_k &=& {1\over{m_k}} \ep\bigl[\tilde
W_{k, m_k}^2 + 2 \tilde W_{k, m_k} (\tilde
W_{k, 2 m_k}-\tilde W_{k, m_k})\bigr]
\nonumber
\\[-8pt]
\\[-8pt]
\nonumber
& = & \sum
_{i=-m_k}^{m_k} \tilde\gamma_{k, i} + 2 \sum
_{i=1}^{m_k} (1-i/m_k) \tilde
\gamma_{k, m_k+i}.
\end{eqnarray}
We now prove that
%
%
\begin{equation}
\label{eqA181047} \nu_k - \sigma^2 = O \Bigl[
\Theta_{m_k, p}+ \min_{l \ge0} \bigl(\Theta_{l, p}
+ l 3^{k (2/p-1)}\bigr) \Bigr],
\end{equation}
which converges to $0$ if $k \to\infty$. Let $\hat X_{k, i} =
T_{3^{k/p}} (X_i)$ and $\hat\gamma_{k,i} = \operatorname{cov}(\hat X_{k,
0},\break
\hat X_{k, i}) = \ep(\hat X_{k, 0} \hat X_{k, i}) - [\ep(\hat
X_{k, 0})]^2$. Note that if $|X_i| \le3^{k/p}$, then $X_i = \hat
X_{k, i}$. Since $X_i \in{\cal L}^p$,
%
%
\begin{eqnarray}
\bigl|\ep(X_0 X_i) - \ep(\hat X_{k, 0} \hat
X_{k, i})\bigr| &=& \bigl|\ep(X_0 X_i \mathbf{
1}_{|X_0| \le3^{k/p}, |X_i| \le3^{k/p}} ) - \ep(\hat X_{k, 0} \hat X_{k, i})
\nonumber\\
&&\hspace*{60pt}{}
+ \ep(X_0 X_i \mathbf{ 1}_{\max(|X_0|, |X_i|) > 3^{k/p}} )\bigr|
\nonumber\\
&\le& \bigl|
\ep(\hat X_{k, 0} \hat X_{k, i} \mathbf{ 1}_{\max(|X_0|, |X_i|) > 3^{k/p}} ) \bigr|
\nonumber
\\[-8pt]
\\[-8pt]
\nonumber
&& {}+ \bigl|\ep(X_0 X_i \mathbf{ 1}_{\max(|X_0|, |X_i|) > 3^{k/p}} )\bigr|
\\
&\le& 2 \ep\bigl[ \bigl(|X_0|+|X_i|\bigr)^2 \mathbf{
1}_{ |X_0| + |X_i| > 3^{k/p}} \bigr]
\nonumber\\
&=& o\bigl(3^{k(2-p)/p}\bigr).\nonumber
\end{eqnarray}
Clearly, we also have $\ep(\hat X_{k, 0}) = o(3^{k(2-p)/p})$. Hence
%
%
\begin{equation}
\label{eqA18651p} \sup_i |\hat\gamma_{k,i} -
\gamma_i| = o\bigl(3^{k(2-p)/p}\bigr).
\end{equation}
For all $j \ge1$, we have $\|W_{k,j} - \tilde W_{k,j} \| \le
j^{1/2} \Theta_{m_k, 2} \le j^{1/2} \Theta_{m_k, p}$. Then
%
%
\begin{equation}
\label{eqA18646p}\qquad \bigl|\ep W_{k,j}^2 - \ep\tilde
W_{k,j}^2\bigr| \le\|W_{k,j} - \tilde
W_{k,j} \| \|W_{k,j} + \tilde W_{k,j} \| \le2 j
\Theta_{m_k, p} \Theta_{0, p}.
\end{equation}
Since $\lim_{j \to\infty} j^{-1} \ep\tilde W_{k,j}^2 =
\sum_{i=-m_k}^{m_k} \tilde\gamma_{k, i}$ and $\lim_{j \to\infty}
j^{-1} \ep W_{k,j}^2 = \sum_{i \in\Z} \hat\gamma_{k,i}$,
(\ref{eqA18646p}) implies that
%
%
\begin{equation}
\label{eqA18647p} \Biggl\llvert \sum_{i=-m_k}^{m_k}
\tilde\gamma_{k, i} - \sum_{i \in\Z} \hat
\gamma_{k,i}\Biggr\rrvert \le2 \Theta_{m_k, p}
\Theta_{0, p}.
\end{equation}
Let the projection operator ${\cal P}_l \cdot= \ep( \cdot| {\cal
F}_l) - \ep( \cdot| {\cal F}_{l-1})$. Then $\hat X_{k,i} = \sum_{l
\in\Z} {\cal P}_l \hat X_{k,i}$. By the orthogonality of ${\cal
P}_l, l \in\Z$, and inequality (\ref{eqA181112}),
%
%
\begin{eqnarray}
\label{eqA181140} |\hat\gamma_{k, i}| &=& \biggl\llvert \sum
_{l \in\Z} \sum_{l' \in\Z} \ep\bigl[({\cal
P}_l \hat X_{k,0}) ({\cal P}_{l'} \hat
X_{k,i} )\bigr] \biggr\rrvert
\nonumber
\\[-8pt]
\\[-8pt]
\nonumber
&\le& \sum
_{l \in\Z} \| {\cal P}_l \hat X_{k,0}\| \| {
\cal P}_{l} \hat X_{k,i} \| \le\sum
_{j=0}^\infty\delta_{j, p}
\delta_{j+i, p}.
\end{eqnarray}
The same inequality also holds for $|\gamma_i|$ and $|\tilde
\gamma_{k, i}|$. For any $0 \le l \le m_k$, we have by~(\ref{eqA181140}) that
%
%
\begin{equation}
\sum_{i=l}^{\infty}\bigl (|\hat
\gamma_{k, i}| + |\tilde\gamma_{k, i}| + |\gamma_i|\bigr)
\le3 \sum_{i=l}^\infty\sum
_{j=0}^\infty \delta_{j, p}
\delta_{j+i, p} \le3 \Theta_{0, p} \Theta_{l,p},
\end{equation}
which entails (\ref{eqA181047}) in view of (\ref{eqA18651p}),
(\ref{eqA18647p}) and (\ref{eqA19513p}).

Recall (\ref{eqA20212p}) and (\ref{eqA161249}) for $\sigma_n^2$.
Now we shall compare $\sigma_n^2$ with
%
%
\begin{equation}
\phi_n = \sum_{k=1}^{h_n-1}
\bigl(3^k - 3^{k-1}\bigr) \nu_k +
\bigl(n-3^{h_n-1}\bigr) \nu_{h_n}.
\end{equation}
Then $\phi_n$ is a piecewise linear function. Observe that, by
(\ref{eqmk1}),
%
%
\begin{equation}
\max_{i \le n} \bigl|\phi_i-\sigma_i^2\bigr|
\le3 \max_{k \le h_n} (m_k \nu_k) = o
\bigl(n^{ (\alpha/p-1)/(\alpha/2 -1)}\bigr).
\end{equation}
By increment properties of Brownian motions, we obtain
%
%
\begin{equation}
\label{eqA20223p} \max_{i \le n} \bigl|\B(\phi_i)-\B\bigl(
\sigma_i^2\bigr)\bigr| = o_\mathrm{ a.s.}\bigl(n^{ (\alpha/p-1)/(\alpha-2)}
\log n\bigr) = o_\mathrm{ a.s.}\bigl(n^{1/p}\bigr).
\end{equation}
Note that by (\ref{eqA181047}), $\phi_i$ is asymptotically linear
with slope $\sigma^2$. Here we emphasize that, under
(\ref{eqsrdsip}), (\ref{eqmk1}), (\ref{eqmap}), a strong
invariance principle with the Brownian motion $\B(\phi_i)$ holds in
view of (\ref{eqA20226}), (\ref{eqA161220}), (\ref{eqA17231p}),
(\ref{eqA161249}), (\ref{eqA20223p}) and Lemma \ref{lemcnstctn}
in the next chapter. However, the approximation $\B(\phi_i)$ is not convenient
for use since~$\phi_i$ is not genuinely linear.

Next, under condition (\ref{eqA18705p}), we shall linearize the
variance function $\phi_i$, so that one can have the readily
applicable form (\ref{eqsipA181103}). Based on the form of
$\phi_i$, we write
%
%
\begin{equation}
\B(\phi_n) = \sum_{k=1}^{h_n-1}
\sum_{j=1}^{3^k - 3^{k-1}} \nu_k^{1/2}
Z_{k,j} + \sum_{j=1}^{n - 3^{h_n-1}}
\nu_{h_n}^{1/2} Z_{h_n,j},
\end{equation}
where $Z_{k, j}$ are i.i.d. standard normal random variables. Define
%
%
\begin{equation}
\B^\ddag(n) = \sum_{k=1}^{h_n-1}
\sum_{j=1}^{3^k - 3^{k-1}} Z_{k,j} + \sum
_{j=1}^{n - 3^{h_n-1}} Z_{h_n,j},
\end{equation}
which is a standard Brownian motion for integer values of $n$. Then we can write
%
%
\begin{equation}
\label{eqA19411p} \B(\phi_n) - \sigma\B^\ddag(n) = \sum
_{i=2}^n b_i Z_i,
\end{equation}
where $(Z_2, Z_3, Z_4, \ldots) = (Z_{1,1}, Z_{1,2}, Z_{2,1},
Z_{2,2}, \ldots, Z_{2, 6}, \ldots, Z_{k,1}, \ldots,\break   Z_{k, 3^k -
3^{k-1}}, \ldots)$ is a lexicographic re-arrangement of $Z_{k,j}$,
and the coefficients $b_n = \nu_{h_n}^{1/2} - \sigma$. Then
%
%
\begin{eqnarray}
\label{eqA18709p} \varsigma_n^2 &=& \bigl\|\B(
\phi_n) - \sigma\B^\ddag(n)\bigr\|^2 = \sum
_{i=2}^n b_i^2
\nonumber
\\[-8pt]
\\[-8pt]
\nonumber
&=& \sum
_{k=1}^{h_n-1} \bigl(3^k -
3^{k-1}\bigr) \bigl(\nu_k^{1/2} - \sigma
\bigr)^2 + \bigl(n - 3^{h_n-1}\bigr) \bigl(
\nu_{h_n}^{1/2} - \sigma\bigr)^2
\end{eqnarray}
and $\varsigma_n^2$ is nondecreasing. If $\lim_{n \to\infty}
\varsigma_n^2 < \infty$, then trivially we have
%
%
\begin{equation}
\label{eqA18745} \B(\phi_n) - \sigma\B^\ddag(n) =
o_\mathrm{ a.s.}\bigl(n^{1/p}\bigr).
\end{equation}
We shall now prove (\ref{eqA18745}) under the assumption that
$\lim_{n \to\infty} \varsigma_n^2 = \infty$. Under the latter
condition, note that we can represent $\B(\phi_n) - \sigma
\B^\ddag(n)$ as another Brownian motion $\B_0(\varsigma_n^2)$, and by
the law of the iterated logarithm for Brownian motion, we have
%
%
\begin{equation}
\label{eqA18749} \mathop{\underline{\overline{\lim}}}_{n \to\infty}
{ {\B(\phi_n) - \sigma\B^\ddag(n)}
\over{\sqrt{2 \varsigma_n^2 \log\log\varsigma_n^2}}} = \pm1 \qquad\mbox{almost surely.}
\end{equation}
Then (\ref{eqA18745}) follows if we can show that
%
%
\begin{equation}
\label{eqA18658p} \varsigma_n^2 \log\log n = o
\bigl(n^{2/p}\bigr).
\end{equation}
Note that (\ref{eqA181047}) and (\ref{eqA18705p}) imply that $3^k
(\nu_k^{1/2} - \sigma)^2 = o(3^{2k/p} / \log k)$, which entails~(\ref{eqA18658p}) in view of (\ref{eqA18709p}).

\section{Some useful lemmas}\label{s4}
\label{secuselem} In this section we shall provide some lemmas that
are used in Section~\ref{secproof}. Lemma \ref{lemcnstctn} is a
``gluing'' lemma, and it concerns how to combine almost sure
convergences in different probability spaces. Lemma \ref{lemasct}
relates truncated and original moments, and Lemma
\ref{lemmomentA14} gives an inequality for moments of maximum sums.

\begin{lemma}
\label{lemcnstctn} Let $(T_{1,n})_{n \ge1}$ and $(U_{1,n})_{n \ge
1}$ be two sequences of random variables defined on the probability
space $(\Omega_1, {\cal A}_1, \pr_1)$ such that $T_{1,n} - U_{1,n}
\to0$ almost surely; let $(T_{2,n})_{n \ge1}$ and $(U_{2,n})_{n
\ge1}$ be another two sequences of random variables defined on the
probability space $(\Omega_2, {\cal A}_2, \pr_2)$ such that $T_{2,n}
- U_{2,n} \to0$ almost surely. Assume that the distributional
equality $(U_{1,n})_{n \ge1} \eqd(T_{2,n})_{n\ge1}$ holds. Then
we can construct a probability space $(\Omega^\dag, {\cal A}^\dag,
\pr^\dag)$ on which we can define $(T_{1,n}')_{n \ge1}$ and
$(U'_{2,n})_{n\ge1}$ such that $(T_{1,n}')_{n \ge1} \eqd
(T_{1,n})_{n \ge1}$, $(U'_{2, n})_{n \ge1} \eqd(U_{2,n})_{n\ge
1}$ and $T_{1,n}' - U'_{2, n} \to0$ almost surely in $(\Omega^\dag,
{\cal A}^\dag, \pr^\dag)$.
\end{lemma}

\begin{pf} Let $\mathbf{ T}_1=(T_{1,n})_{n\ge1}$, $\mathbf{
U}_1=(U_{1,n})_{n\ge1}$, $\mathbf{ T}_2=(T_{2, n})_{n\ge1}$,
$\mathbf{
U}_2=\break  (U_{2,n})_{n\ge1}$; let $\mu_{\mathbf{ T}_1 |\mathbf{ U}_1}$ and $\mu
_{\mathbf{
U}_2 | \mathbf{T}_2}$ denote, respectively, the conditional distribution of $ \mathbf{T}_1 $ given
$ \mathbf{U}_1$ and the conditional distribution of $ \mathbf{U}_2 $ given
$ \mathbf{T}_2$. Let $(\Omega^\dag, {\cal F}^\dag, P^\dag)$ be a
probability space on which there exists a vector $ \mathbf{U}_1'$
distributed as $ \mathbf{U}_1$. By enlarging $(\Omega^\dag, {\cal F}^\dag,
P^\dag)$ if necessary, there exist random vectors $ \mathbf{T}_1'$ and
$ \mathbf{U}_2'$ on this probability space such that the conditional
distribution of $ \mathbf{T}_1'$ given $ \mathbf{U}_1'$ equals $\mu_{\mathbf{ \mathbf{T}_1 |
\mathbf{U}_1}}$, and the conditional distribution of $ \mathbf{U}_2'$ given $
\mathbf{U}_1'$ equals $\mu_{\mathbf{ U}_2 | \mathbf{T}_2}$. Then by $\mathbf{ U}_1 \eqd
\mathbf{
T}_2$ we have $ (\mathbf{T}_1', \mathbf{U}_1') \eqd (\mathbf{T}_1, \mathbf{U}_1)$ and
$(\mathbf{U}_1', \mathbf{U}_2') \eqd (\mathbf{T}_2, \mathbf{U}_2)$, so that for the components we
have $T_{1,n}'-U_{1,n}'\to0$ a.s. and $U_{1,n}'-U_{2,n}'\to0$
a.s., so that $T_{1,n}'-U_{2,n}' \to0$ a.s.
\end{pf}

\begin{lemma}
\label{lemasct} Let $X \in{\cal L}^p$, $2 < p < \alpha$. Then
there exists a constant $c = c_{\alpha, p}$ such that
%
%
\begin{equation}
\label{eqCT43}\qquad \sum_{i=1}^\infty3^i
\pr\bigl(|X| \ge3^{i/p}\bigr) + \sum_{i=1}^\infty3^i
\ep\min\bigl(\bigl|X/3^{i/p}\bigr|^\alpha, \bigl|X/3^{i/p}\bigr|^2
\bigr) \le c \ep\bigl(|X|^p\bigr).
\end{equation}
\end{lemma}
\begin{pf} That the first sum is finite follows from
%
%
\begin{equation}
\label{eqCT431} \sum_{i=1}^\infty3^i
\pr\bigl(|X| \ge3^{i/p}\bigr) \le3 \sum_{i=1}^\infty
\int_{3^{i-1}}^{3^i} \pr\bigl(|X|^p > u\bigr)
\,d u \le3 \ep\bigl(|X|^p\bigr).
\end{equation}
For the second one, let $q_i = \pr(3^{i-1} \le|X|^p < 3^i)$. Then
%
%
\begin{eqnarray}
\label{eqCT432} \sum_{i=1}^\infty3^i
\ep\bigl(\bigl|X/3^{i/p}\bigr|^2 \mathbf{ 1}_{|X|^p\ge3^i}\bigr) &
\le& \sum_{i=1}^\infty3^i \sum
_{j=1+i}^\infty 3^{(j-i)2/p}
q_j
\nonumber\\
&=& \sum_{j=2}^\infty
\sum_{i=1}^{j-1} 3^i
3^{(j-i)2/p} q_j
\\
&=& c_1 \sum
_{j=2}^\infty3^{j} q_j \le
c_1 \ep\bigl(|X|^p\bigr)\nonumber
\end{eqnarray}
for some constant $c_1$ only depending on $p$ and $\alpha$.
Similarly, there exists $c_2$ such that
\begin{eqnarray*}
\sum_{i=1}^\infty3^i \ep
\bigl(\bigl|X/3^{i/p}\bigr|^\alpha\mathbf{ 1}_{|X|^p < 3^i}\bigr) &\le&
\sum_{i=1}^\infty3^i \sum
_{j=-\infty}^i 3^{(j-i)\alpha/p} q_j
\\
&=&
\sum_{j=-\infty}^\infty\sum
_{i=\max(1,j)}^\infty 3^{i(1-\alpha/p)} 3^{j \alpha/p}
q_j \le c_2 \ep\bigl(|X|^p\bigr).
\end{eqnarray*}
For the last relation, we consider the two cases $\sum_{j=-\infty}^0$
and $\sum_{j=1}^\infty$ separately. The lemma then follows from
(\ref{eqCT431}) and (\ref{eqCT432}). It is easily seen that
(\ref{eqCT43}) also holds with the factor $3$ therein replaced by
any $\theta> 1$. In this case the constant $c$ depends on $p,
\alpha$ and $\theta$.
\end{pf}

\begin{lemma}
\label{lemmomentA14} Recall (\ref{eqsrdsip}) and (\ref{eqmk1})
for $\Xi_{\alpha, p}$ and $M_{\alpha, p}$, respectively, and
(\ref{eqA14725p}) for $W_{k, l}$. Then there exists a constant $c$,
only depending of $\alpha$ and $p$, such that
%
%
\begin{equation}
\label{eqA14728p} \sum_{k=1}^\infty
{ {3^k}\over{m_k}} { {\ep(\max_{1\le l \le m_k} |W_{k, l}|^\alpha)}
\over{ 3^{k\alpha/p} }} \le c M_{\alpha, p}
\Theta_{0, 2}^\alpha + c \Xi_{\alpha, p}^\alpha + c
\|X_1\|_p^p.
\end{equation}
\end{lemma}

\begin{pf}Recall
(\ref{eqtruncSA14738}) for the functional dependence measure
$\delta_{k, j, \iota}$. Since $T_a$ has Lipschitz constant $1$, we
have
%
%
\begin{eqnarray}
\label{eqtruncSA14748} \delta_{k, j, \iota}^\iota &\le& \ep\bigl[\min
\bigl(2 \times3^{k/p}, |X_i-X_{i, \{i-j\}}|
\bigr)^\iota\bigr]
\nonumber
\\[-8pt]
\\[-8pt]
\nonumber
&\le& 2^\iota\ep\bigl[\min
\bigl(3^{k/p}, |X_j-X_{j, \{0\}}|\bigr)^\iota
\bigr].
\end{eqnarray}
We shall apply the Rosenthal-type inequality in \citet{LiuHanWu}: there exists a constant $c$, only depending on $\alpha$,
such that
%
%
\begin{eqnarray}
\label{eqA14808} \Bigl\llVert \max_{1\le l \le m_k} |W_{k, l}|
\Bigr\rrVert _{\alpha} &\le& c m_k^{1/2} \Biggl[ \sum
_{j=1}^{m_k} \delta_{k, j, 2} + \sum
_{j=1+m_k}^\infty\delta_{k, j, \alpha} + \bigl\|
T_{3^{k/p}} (X_1)\bigr\|_2 \Biggr]
\nonumber\\
& &{}+ c
m_k^{1/\alpha} \Biggl[ \sum_{j=1}^{m_k}
j^{1/2-1/\alpha} \delta_{k, j, \alpha} + \bigl\| T_{3^{k/p}} (X_1)
\bigr\|_\alpha \Biggr]
\\
&\le& c(I_k + \mathit{I I}_k + \mathit{I I
I}_k),\nonumber
\end{eqnarray}
where
%
%
\begin{eqnarray}
\label{eqA14826} I_k &=& m_k^{1/2} \sum
_{j=1}^\infty\delta_{j, 2} +
m_k^{1/2} \|X_1\|_2,
\nonumber\\
\mathit{I
I}_k &=& m_k^{1/\alpha} \sum
_{j=1}^\infty j^{1/2-1/\alpha} \delta_{k, j, \alpha},
\\
\mathit{I I I}_k &=& m_k^{1/\alpha} \bigl\| T_{3^{k/p}}
(X_1)\bigr\|_\alpha.\nonumber
\end{eqnarray}
Here we have applied the inequality $\delta_{k, j, 2} \le\delta_{j,
2}$, since $T_a$ has Lipschitz constant $1$. Since
$\sum_{j=1}^\infty\delta_{j, 2} + \|X_1\|_2 \le2 \Theta_{0, 2}$,
by (\ref{eqmk1}), we obtain the upper bound $c M_{\alpha, p}
\Theta_{0, 2}^\alpha$ in (\ref{eqA14728p}), which corresponds to
the first term $I_k$ in (\ref{eqA14808}). For the third term
$\mathit{I I I}_k$, we obtain the bound $c \|X_1\|_p^p$ in
(\ref{eqA14808}) in view of Lemma \ref{lemasct} by noting that
$|T_{3^{k/p}} (X_1)| \le\min(3^{k/p}, |X_1|)$ and
$\min(|v|^\alpha, v^2) \ge\min(|v|^\alpha, 1)$.

We shall now deal with $\mathit{I I}_k$. Let $\beta= \alpha/(\alpha-1)$, so
that $\beta^{-1} + \alpha^{-1} = 1$; let $\lambda_j = (j^{1/2 -
1/\alpha} \delta_{j, p}^{p/\alpha})^{-1/\beta}$. Recall
(\ref{eqsrdsip}) for $\Xi_{\alpha, p}$. By H\"older's inequality,
%
%
\begin{equation}
\Biggl(\sum_{j=1}^\infty j^{1/2-1/\alpha}
\delta_{k, j, \alpha} \Biggr)^\alpha \le\Xi_{\alpha, p}^{\alpha/\beta}
\sum_{j=1}^\infty\lambda_j^\alpha
\bigl(j^{1/2-1/\alpha} \delta_{k, j,
\alpha}\bigr)^\alpha.
\end{equation}
Hence, by (\ref{eqtruncSA14748}) and Lemma \ref{lemasct}, we
complete the proof of (\ref{eqA14728p}) in view of
%
%
\begin{eqnarray}
\label{eqA14903} \sum_{k=1}^\infty
{ {3^k}\over{m_k}} { {\mathit{I I}_k^\alpha} \over{3^{\alpha k/p}}} &\le& \sum
_{k=1}^\infty3^{k-k\alpha/p} \Xi_{\alpha, p}^{\alpha/\beta}
\sum_{j=1}^\infty\lambda_j^\alpha
\bigl(j^{1/2-1/\alpha} \delta_{k, j,
\alpha}\bigr)^\alpha
\nonumber\\
&=&
\Xi_{\alpha, p}^{\alpha/\beta} \sum_{j=1}^\infty
\lambda_j^\alpha j^{\alpha/2-1} \sum
_{k=1}^\infty3^{k-k\alpha/p} \delta_{k, j, \alpha}^\alpha
\\
&\le& \Xi_{\alpha, p}^{\alpha/\beta} \sum_{j=1}^\infty
\lambda_j^\alpha j^{\alpha/2-1} c_{\alpha, p}
\delta_{j, p}^p = c_{\alpha, p} \Xi_{\alpha, p}^{\alpha}.\nonumber
\end{eqnarray}
\upqed\end{pf}
\section*{Acknowledgment} We thank the anonymous referee for
his/her helpful comments that have improved the paper.
We also thank F. Merlev\`{e}de and E. Rio for
pointing out an error in an earlier version of the paper.

%



\printaddresses

\end{document}